\documentclass[a4paper]{article}
\usepackage{amssymb}
\usepackage{amsmath}
\usepackage{amsthm}
\usepackage{color}
\usepackage{geometry}
\usepackage{graphicx}
\usepackage{epsf}
\usepackage{subfigure}
\usepackage{multirow}
\usepackage{mathrsfs}
\newtheorem{thm}{Theorem}[section]
\newtheorem{lam}[thm]{Lemma}

\newtheorem{Rem}[thm]{Remark}


\title{Normalized solutions for nonlinear Schr\"odinger equations on graphs}

\begin{document}

\date{2023.1.16}

\author{Yunyan Yang$^{1}$\ \ \ Liang Zhao$^{2*}$\\
\it\small  $^{1}$  School of Mathematics, Renmin University of China, Beijing 100872, China\\
\it\small  $^{2}$  School of Mathematical Sciences, Key Laboratory of Mathematics and Complex Systems of MOE, \\
\it\small Beijing Normal University, Beijing 100875, China\\
\it\small $^*$ Corresponding author, E-mail: liangzhao@bnu.edu.cn
}

\maketitle

\begin{abstract}
We are concerned with the nonlinear Schr\"odinger equation with an $L^2$ mass constraint on both finite and locally finite graphs and prove that the equation has a normalized solution by employing variational methods. We also pay attention to the behaviours of the normalized solution as the mass constraint tends to $0^+$ or $+\infty$ and give clear descriptions of the limit equations. Finally, we provide some numerical experiments on a finite graph to illustrate our theoretical results.



\end{abstract}

\section{Introduction}

Graph theory and PDE are two research fields with fundamental theoretical significance and application value. As an intersection of them, partial differential equations on graphs have become a very attractive topic these years. For example, Chow, Li and Zhou considered the entropy dissipation of Fokker-Planck equation \cite{Chow1} and optimal transport \cite{Chow2} on graphs. Horn, Lin, Liu and Yau \cite{Horn} proved Li–Yau-type estimates for solutions of the heat equation on graphs. In particular, Grigor'yan, Lin and Yang \cite{Gri1,Gri2,Gri3} first studied the variational structure of a series of nonlinear elliptic equations on graphs and pointed out that the required Sobolev spaces are pre-compact, which made it possible to apply variational methods to the existence of solutions. Along this line, there are many follow-up research progresses. The existence and limit behaviour of nontrivial solutions to Schr\"odinger equation and system with potential well were obtained by Zhang, Zhao, Han, Shao and Xu \cite{ZhangZhao,Han1,Han2,XuZhao}. Similar problems on infinite metric graphs were studied by Akduman and Pankov \cite{Akduman-Pankov}. The results of Grigor'yan, Lin and Yang on Kazdan-Warner equation was extended by Ge and Jiang \cite{GeJiang} to infinite graphs and by Keller and Schwarz \cite{KellerSchwarz} to canonically compactifiable graphs. Applications of the degree theory on Kazdan-Warner equation were also discussed by Huang, Wang, Yang, Sun and Liu \cite{HuangWangYang,SunWang}. The semilinear heat equation on locally finite graphs were studied by Huang, Lin and Wu \cite{Huang, Lin1}. Multiplicity of solutions for nonlinear equations on finite graphs were investigated by Liu, Yang, Chao and Hou \cite{LiuYang,ChaoHou}. For other related works, readers are referred to \cite{ZhangLin1,ZhangLin2,LinYang1,LinYang2} and the references therein.
 
In this note, we are concerned with the following nonlinear equation
\begin{equation}\label{shrodinger}
-\Delta u+\lambda u=f\left (u\right ),	
\end{equation}
where $\lambda\in \mathbb{R}$ and $f$ is a nonlinearity satisfying some necessary conditions. When the unknown real-valued function $u$ is defined on an Euclidean space $\mathbb{R}^n$, this equation is motivated by searching for standing wave solutions of the Schr\"odinger or Klein-Gordon equations. Indeed, for the Schr\"odinger equation
\begin{equation}\label{cshrodinger}
	i \Phi_t+\Delta \Phi+g(|\Phi|^2)\Phi=0,	
\end{equation}
where $\Phi(x,t):\mathbb{R}^n\times\mathbb{R}^+\rightarrow \mathbb{C}$, when we search for its standing wave solutions, that is, $\Phi(x,t)=e^{i \lambda t}u(x)$, it leads to solving equation \eqref{shrodinger} with $f(u)=g(|u|^2)u$. For this reason, the equations \eqref{shrodinger} is related to several important physical models, such as the deep-water waves and the quantum states of particles in nonlinear optics and Bose-Einstein condensates. We usually call the equation \eqref{shrodinger} a Schr\"odinger type equation and it has been extensively studied.

When we consider the equation \eqref{shrodinger} with a given constant $\lambda\in \mathbb{R}$, one can apply variational or some other topological methods to obtain kinds of properties of the equation, such as the existence, multiplicity and asymptotic behaviour, etc. There are too many related literatures and we can not summarize all of them here. The readers can refer to \cite{Lions1,Lions2,Brezis,Cao,Rab,Alves,Chen} and the references therein. On the other hand, it is obvious that the $L_2$ norm $\Vert \Phi\Vert_2$ of the standing wave solution $\Phi(x,t)=e^{i\lambda t}u(x)$ equals to $\Vert u\Vert_2$, where $u(x)$ is a solution of \eqref{shrodinger}. From the physical point of view, $\Vert \Phi\Vert_2$ represents the mass of particles or the normalization constant of a probability distribution function. Therefore, the solution of \eqref{shrodinger} with a prescribed $L_2$ norm is particularly interesting because it corresponds to a standing wave solution of \eqref{cshrodinger} which preserves its mass along the time evolution. We call
\begin{equation}\label{mass}
	\left \{
	\begin{array}{cl}
		-\Delta u+\lambda u=f\left (u\right ) \ \ \text{in}\ \  \mathbb{R}^n\\
		\Vert u\Vert_2^2=m
	\end{array}
	\right .
\end{equation}
a mass preserved problem and its solution a normalized solution. We point out that $\lambda\in \mathbb{R}$ in \eqref{mass} arises as a Lagrange multiplier instead of a given constant in \eqref{shrodinger}. In the last few years, the mass preserved problem on $\mathbb{R}^n$ has been a hot topic and the readers can refer to \cite{Jean3,Shi,Bon,Hira,Iko,Ste1,Jean2,Soa,BarMol,Pel} and the references therein for more details. Based on the above backgrounds, we aim to study the discrete mass preserved problem on a graph $G=\left(V,E\right)$. In this situation, the normalized solution on a graph can be regarded as a discrete standing wave solution or a particle's probability distribution function in a discrete space. As far as the authors know, this is the first research work about the mass preserved problem on graphs.

To describe our problem in details, we first introduce some basic concepts and assumptions. Let $G=\left(V,E\right)$ be a connected graph, where $V$ denotes the set of vertices and $E$ denotes the set of edges. A graph $G$ is finite if the number of vertices and edges are both finite. While $G$ is locally finite if for any $x\in V$, there are only finite $y\in V$ such that $xy\in E$. We always assume that the weight $w_{xy}\in \mathbb{R}$ of any edge $xy \in E$ is positive and satisfies $w_{xy}=w_{yx}$. The degree of a vertex $x\in V$ is defined as $deg(x):=\sum_{xy\in E}w_{xy}$. The measure $\mu :V\rightarrow \mathbb{R}^{+}$ is a finite positive function on $G$ and the $\mu$-Laplacian of any function $u: V\rightarrow \mathbb{R}$ at $x\in V$ is
\begin{displaymath}
\Delta u( x ) :=\frac{1}{\mu ( x ) } \sum_{y\sim x}^{} w_{xy}( u(y)-u(x)  ),
\end{displaymath}
where $y\sim x$ means that there exists an edge $xy\in E$ connecting $x$ and $y$. The gradient form at $x\in V$ of two functions $u$ and $v$ is
\begin{displaymath}
\Gamma (u,v)(x):=\frac{1}{2\mu(x)}\sum_{y\sim x}^{}  w_{xy}(u(y)-u(x))(v(y)-v(x)).
\end{displaymath}
In particular, we use $\Gamma(u)$ to denote $\Gamma(u,u)$ and the length of the gradient for $u$ at $x\in V$ is
\begin{displaymath}
\left | \nabla u \right | (x):=\sqrt{\Gamma(u)(x)}=\left(\frac{1}{2\mu(x)}\sum_{y\sim x}^{}w_{xy} \big (u(y)-u(x)\big )^{2} \right )^{1/2}.
\end{displaymath}
To do calculus of variations on graphs, we also need the concepts of integral and function spaces. The integral of a function $u$ over $V$ is
\begin{displaymath}
\int_{V}^{}ud\mu=\sum_{x\in V}^{}  \mu(x)u(x).
\end{displaymath}
For any $q>0$, $L_q(V)$ is a linear space of functions $u:V\rightarrow \mathbb{R}$ with the norm
$$\Vert u\Vert_q=\left(\int_{V}|u|^q d\mu\right)^{1/q}<\infty.$$
$L_\infty(V)$ is the space with the norm
$$\Vert u\Vert_\infty=\sup_V |u|<\infty.$$
Moreover, $H(v)$ denotes the space with the norm
\begin{displaymath}
\Vert u \Vert _{H}= \left ( \int_{V}^{} \left ( \left | \nabla u \right |^{2}+u^{2}  \right ) d\mu  \right )^{1/2}<\infty.
\end{displaymath}
Obviously, $H(v)$ is a Hilbert space with the inner product
$$\left\langle  u,v\right\rangle=\int_{V}\left(\Gamma(u,v)+uv\right) d\mu,$$
for any $u,v\in H(V)$.

We deal with the mass preserved problem on both finite and locally finite graphs by using different methods. Now, let us present our main results for these two cases respectively. On a connected finite graph $G=(V,E)$, we consider the following nonlinear Schr\"odinger equation
\begin{equation}\label{eqf}\left\{\begin{array}{lll}
		-\Delta u+\lambda_m u=u^{p-1}&{\text in}& V\\
		u>0&{\text in}& V\\
		\int_Vu^2d\mu=m
	\end{array}\right.
\end{equation}
where $p>2$, $m>0$ is the mass constraint and $\lambda_m$ is the Lagrange multiplier corresponding to $m$. For the existence of solutions to this equation, we have

\begin{thm}\label{existf}
For any $m>0$ and $p>2$, the equation \eqref{eqf} on a connected finite graph $G=(V,E)$ has a positive solution $u_m$.
\end{thm}

It is also interesting to explore the behaviour of the solution $u_m$ as the evolution of the mass. Suppose that $u_m$ is the solution of \eqref{eqf} corresponding to mass $m$. We have
\begin{thm}\label{limitf}
(i) If $m\rightarrow 0^+$, up to a subsequence, we have $u_m=\sqrt{m}v+o(\sqrt{m})$, where either $v\equiv |V|^{-\frac{1}{2}}$ on $V$, or $v:V\rightarrow\mathbb{R}$  satisfies
\begin{equation*}\label{m0f}
	\left\{\begin{array}{lll}
		-\Delta v=\lambda_0 v&{\text in}& V\\
		v>0&{\text in}& V\\
		\int_Vv^2d\mu=1
	\end{array}\right.
\end{equation*}
for some constant $\lambda_0>0$. Here $|V|=\int_V 1d\mu$ denotes the volume of $V$.

(ii) If $m\rightarrow+\infty$, up to a subsequence, we have $u_m=\sqrt{m}(w+o_m(1))$, where $w:V\rightarrow\mathbb{R}$ satisfies
$$w^{p-1}-\Vert w\Vert_{p}^p w=0.$$
\end{thm}

If $G=(V,E)$ is a locally finite graph, we consider the following nonlinear Schr\"odinger equation
\begin{equation}\label{eql}
	\left\{\begin{array}{lll}
		-\Delta u+\lambda_m hu=u^{p-1} &{\text in}& V\\
		u>0&{\text in}& V\\
		\int_Vhu^2d\mu=m
	\end{array}\right.
\end{equation}
Since a locally finite graph can contain infinite vertices and edges, the problem is more complicated and we need some additional assumptions as follows.

(c1) The positive measure $\mu:V\rightarrow \mathbb{R}$ has a uniformly lower bound, i.e., there exists some constant $\mu _{\min}>0$ such that $\mu \left ( x \right ) \geqslant \mu _{\min}$ for any $x\in V$.

(c2) The function $h:V\rightarrow \mathbb{R}$ satisfies $\inf_{x\in V}h(x)=h_0>0$ and $h(x)\rightarrow+\infty$ as $\rho(x)\rightarrow\infty$, where $\rho(x)=d(x,O)$ denotes the distance between any vertex $x\in V$ and a given vertex $O\in V$. 

(c3) For the function $h$, during the proof of the existence of solutions to \eqref{eql}, we also assume that
\begin{equation*}
	h(O)<m\left\{\frac{2}{p}\frac{\mu(O)^{2-\frac{p}{2}}}{{\text deg}(O)}\right\}^{\frac{2}{p-2}}.
\end{equation*}

\begin{thm}\label{existl}
If $G=(V,E)$ is a connected locally finite graph satisfying (c1) and $h$ is a function satisfying (c2) and (c3), the equation \eqref{eql} has a positive solution $u_m$ for any $m>0$ and $p>2$.
\end{thm}

For the behaviour of $u_m$ on a locally finite graph as the mass $m$ varying, we have

\begin{thm}\label{limitl}
Suppose that (c1)-(c3) are satisfied, $u_m$ is a solution of \eqref{eql} with a mass constraint $m$ and $\lambda_m\geq-c$ for some constant $c\geq 0$.

(i) If $m\rightarrow 0+$, up to a subsequence, $v_m=u_m/{\sqrt{m}}$ converges to a function $v$ uniformly on $V$, where either $v\equiv 0$, or $v$ is a solution of the equation
\begin{equation*}
	\left\{\begin{array}{lll}
	-\Delta v=\lambda_0 hv&{\text in}& V\\
	0<\int_Vhv^2d\mu\leq 1
    \end{array}\right.
\end{equation*}
for some constant $\lambda_0>0$. In particular, if $c=0$, we have $v\equiv 0$ on $V$.
	
(ii) If $m\rightarrow+\infty$, up to a subsequence, $w_m=u_m/\sqrt{m}$ converges to a function $w$ uniformly on $V$, where either $w\equiv 0$, or $w$ is a solution of the equation
$$|w|^{p-2}w=\lambda_\infty hw,$$
where
\begin{equation}\label{lambdainfty}
	\lambda_\infty=\frac{\int_V|w|^pd\mu}{\int_Vhw^2d\mu}.
\end{equation}
\end{thm}

\begin{Rem}
For several equations that seem to be similar to the equations we are discussing, there are some existence results on graphs. Zhang \cite{ZhangLin1} proved that the equation
$$-\Delta u+hu=\lambda f(u)$$
has a solution on graphs. Stefanov, Ross and Kevrekidis \cite{Ste2} obtained the existence for the equation
\begin{equation*}
	\left\{\begin{array}{lll}
		-\Delta u+au^{p-1}=\lambda u^{q-1}\\
		u>0\\
		\|u\|_q^q=m
	\end{array}\right.
\end{equation*}
on a lattice graph. Very recently, Hua, Li and Wang \cite{Hua} extended this result and proved that on a Cayley
graph of a discrete group of polynomial growth with the homogeneous dimension $N\geq 2$, the equation has a solution for any $2\leq p<q<+\infty$, $a>0$ and $m>0$. But these equations are all different from our equations \eqref{eqf} and \eqref{eql} and their methods are also very different from ours.
\end{Rem}

The paper is organized as follows. We first deal with the finite graph case in Section 2. Theorem \ref{existf} and \ref{limitf} are proved in this section. Section 3 is devoted to the locally finite graph case and we prove Theorem \ref{existl} and \ref{limitl} in this section. In Section 4, we give some numerical experiments on a finite graph to illustrate our theorems.

\section{The finite graph case}

In this section, we always assume that $G=(V,E)$ is a connected finite graph. Some fundamental tools for calculus of variations on graphs, such as formulas of integration by parts, the definition of a weak solution, can be found in many previous literatures, such as \cite{Han1} and \cite{ZhangZhao}, and we omit the details here.

The functional corresponding to the equation \eqref{eqf} is
\begin{equation}\label{funcf}
J(u)=\frac{1}{2}\int_V|\nabla u|^2d\mu-\int_V|u|^pd\mu,\quad\forall u\in H(V).
\end{equation}
For the finite graph case, the existence of solutions for \eqref{eqf} can be proved by searching for a global minimizer of the functional \eqref{funcf} in
$$
S_m:=\left\{u\in H(V): \|u\|_{L^2(V)}^2=m\right\}.
$$

\vspace{0.5cm}

\noindent {\bf Proof of Theorem \ref{existf}.}
Suppose that there are $k\in \mathbb{Z}^+$ vertices in $V$, i.e., $V=\{x_1,\cdots,x_k\}$. We use the notations $y_i=u(x_i)$, $i=1,\cdots, k$, and $\mathbf{y}=(y_1,\cdots,y_k)\in\mathbb{R}^k$ for any function $u\in H(V)$. Since there are only finite vertices in $V$, $H(V)$ is identified with the finite dimensional vector space
$\mathbb{R}^k$. In particular, $u\in S_m$ if and only if
$$\mathbf{y}\in A_m=\left\{y\in\mathbb{R}^k:\sum_{i=1}^k\mu(x_i)y_i^2=m\right\}.$$
Define a function $F:\mathbb{R}^k\rightarrow\mathbb{R}$ by
$$
F(\mathbf{y})=J(u)=\frac{1}{4}\sum_{x\in V,z\sim x}w_{zx}(u(z)-u(x))^2-\sum_{x\in V}\mu(x)|u|^p(x).
$$
Clearly $F\in C^1(\mathbb{R}^k,\mathbb{R})$. Observing
$$\inf_{u\in S_m}J(u)=\inf_{y\in A_m}F(\mathbf{y})$$
and noting that $A_m$ is a compact subset of $\mathbb{R}^k$, we certainly can find some $\mathbf{y}_0=(y_{0_1},\cdots,y_{0_k})\in A_m$ such that
$F(\mathbf{y}_0)$ achieves the above infimum. Define a function $u_0:V\rightarrow\mathbb{R}$ by $u_0(x_i)=y_{0_i}$, $i=1,\cdots,k$. Obviously, there hold $u_0\in S_m$ and
$$J(u_0)=\inf_{u\in S_m}J(u).$$

Moreover, since $|\nabla |u||^2\leq |\nabla u|^2$, we have
$$J(|u|)\leq J(u),\quad\forall u\in H(V).$$
Consequently, there hold $|u_0|\in S_m$ and
$$J(|u_0|)=\inf_{u\in S_m}J(u).$$
With no loss of generality, we can assume $u_0\geq 0$. By a straightforward calculation, the Euler-Lagrange equation of $u_0$ is given as
\begin{equation}\label{eqfpos}
	\left\{\begin{array}{lll}
	-\Delta u_0+\lambda_m u_0=u_0^{p-1}&{\text in}& V\\
	u_0\geq 0&{\text in}& V\\
	\int_Vu_0^2d\mu=m,
    \end{array}\right.
\end{equation}
where $\lambda_m\in\mathbb{R}$ is the Lagrange multiplier defined by
$$\lambda_m=\frac{1}{m}\left\{\int_V|u_0|^p d\mu-\int_V|\nabla u_0|^2d\mu\right\}.$$
We claim that $u_0(x)$ is strictly positive for any $x\in V$. Indeed, suppose that there exists some
$x_0\in V$ such that
$u_0(x_0)=0=\min_{x\in V} u_0(x)$. If there exists some $z\sim x_0$ such that $u(z)>0$, we have
$$0>\frac{1}{\mu(x_0)}w_{zx_0}(u(x_0)-u(z))\geq -\Delta u_0(x_0)=u_0^{p-1}(x_0)-\lambda_m u_0(x_0)= 0,$$
which is impossible. Hence we have $u_0(z)=0$ if $z\sim x_0$. Since $G$ is connected, by repeating this process for finite times, we can conclude that $u_0\equiv 0$ on $V$, which contradicts the mass constraint $\|u_0\|_{L^2(V)}^2=m>0$.

By the above discussions, we get a strictly positive solution of \eqref{eqf} which is a global minimizer of the functional \eqref{funcf} and the theorem is proved.
$\hfill\Box$

\begin{Rem}
We can also prove Theorem \ref{existf} for some more general nonlinearity $f(x,t): V\times\mathbb{R}\rightarrow\mathbb{R}$. For example, if $f(x,t)$ is continuous in $t\in\mathbb{R}$, $|f(x,t)|\leq f(x,|t|)$ and $f(x,t)\not\equiv 0$
for $(x,t)\in V\times[0,+\infty)$, the above proof can be carried on with only minor modifications.
\end{Rem}

Next, we consider the limit behaviour of the solution to \eqref{eqf} with the evolution of the mass $m$.
\vspace{0.5cm}

\noindent{\bf Proof of Theorem \ref{limitf}.}
Suppose that $u_m$ is a solution of \eqref{eqf} with the mass constraint $m>0$. Obviously for $v_m=u_m/\sqrt{m}$, we have
$$\int_Vv_m^2d\mu=1.$$

{\bf Case 1. $m\rightarrow 0^+$.}
Since $H(V)$ is finite dimensional and pre-compact, there exists a subsequence, which is still denoted by $\{v_m\}$, such that $v_m(x)\rightarrow v(x)$ for all $x\in V$ as $m\rightarrow 0^+$. Direct computations give us that $v_m$ satisfies the equation
\begin{equation}\label{eqfmv}
	-\Delta v_m+\lambda_mv_m={m}^{\frac{p}{2}-1}v_m^{p-1}.
\end{equation}
Since $p>2$, we have
\begin{equation}\label{r0}
	{m}^{\frac{p}{2}-1}v_m^{p-1}(x)\rightarrow 0\,\,{\rm uniformly \   for\ } x\in V,{\rm \  as\ } m\rightarrow 0^+.
\end{equation}
In view of (\ref{eqfmv}), we have
\begin{eqnarray}\nonumber
\lambda_m&=&m^{\frac{p}{2}-1}\int_Vv_m^pd\mu-\int_V|\nabla v_m|^2d\mu\\
&=&-\int_V|\nabla v|^2d\mu+o_m(1),\label{r01}
\end{eqnarray}
where $o_m(1)\rightarrow 0$ as $m\rightarrow 0^+$. Combining (\ref{eqfmv})-(\ref{r01}), we conclude that
$$\left\{\begin{array}{lll}
	-\Delta v=\lambda_0 v&{\rm in}& V\\
	v\geq 0&{\rm in}& V\\
	\int_Vv^2d\mu=1,
\end{array}\right.
$$
where
$$\lambda_0=\int_V|\nabla v|^2d\mu=-\lim_{m\rightarrow 0^+}\lambda_m.$$
If $\lambda_0=0$, $v$ shall be a harmonic function and the maximum principle implies that $v\equiv c$ for some constant $c\in \mathbb{R}$. Consequently, $\int_Vv^2d\mu=1$ gives us that $v(x)\equiv |V|^{-\frac{1}{2}}$ for all $x\in V$.
While for the case $\lambda_0\not=0$, the maximum principle implies that $v(x)>0$ for any $x\in V$.
Moreover, $\lambda_0=\int_V|\nabla v|^2d\mu$ gives us that $\lambda_0>0$.
\vspace{0.5cm}

{\bf Case 2. $m\rightarrow +\infty$.}
In this case, \eqref{eqfmv} still holds for $w_m=u_m/\sqrt{m}$. Without loss of generality, we assume that $w_m$ converges to some function $w$ uniformly in $V$ as $m\rightarrow+\infty$. In this case, instead of (\ref{r01}), we have
$$
\lambda_m=m^{\frac{p}{2}-1}(\|w\|_{p}^p+o_m(1)).
$$
In view of (\ref{eqfmv}), we obtain
$$
-\Delta v_m=m^{\frac{p}{2}-1}\left(w^{p-1}-\|w\|_{p}^p w+o_m(1)\right),
$$
which immediately leads to
$$
w^{p-1}-\Vert w\Vert_{p}^p w=0.$$
	This ends the proof of the theorem.
$\hfill\Box$

\section{The locally finite graph case}

In this section, we always suppose that $G=(V,E)$ is a connected locally finite graph. $C_c(V)$ denotes the set of functions $u:V\rightarrow \mathbb{R}$ such that $\left\{ x\in V:u\left ( x \right )\ne 0 \right\}$ is of finite cardinality. Let $H_0(V)$ be the completion of $C_c(V)$ under the norm
$$\|u\|_{H_0}=\left(\int_V(|\nabla u|^2+u^2)d\mu\right)^{1/2}$$
and define the Hilbert space $\mathcal{H}$ as
$$
\mathcal{H}=\left\{u\in H_0(V):\int_V(|\nabla u|^2+hu^2)d\mu<+\infty\right\}.
$$
We use
$$
\mathcal{J}(u)=\frac{1}{p}\int_V|u|^pd\mu-\frac{1}{2}\int_V|\nabla u|^2d\mu, \ \ \forall u\in \mathcal{H}.
$$
as the functional corresponding to \eqref{eql} and will solve the equation \eqref{eql} by considering the maximizing problem
$$
\Lambda_m=\sup\left\{\mathcal{J}(u):u\in \mathcal{B}_m\right\},
$$
where
$$
\mathcal{B}_m=\left\{u\in\mathcal{H}:\int_Vhu^2d\mu\leq m\right\}.
$$
We also use $\mathcal{S}_m$ to denote the following set of functions
$$\mathcal{S}_m=\{u\in\mathcal{H}: \int_Vhu^2d\mu=1\}\subset \mathcal{B}_m.$$
Before the proof of Theorem \ref{existl}, we first present a useful lemma for the lower bound of $\Lambda_1$.

\begin{lam}\label{bound}
	If the conditions (c1)-(c3) are satisfied, we have $\Lambda_1=\sup\left\{\mathcal{J}(u):u\in \mathcal{B}_1\right\}>0$.
\end{lam}

\noindent {\bf Proof.}
To prove the lemma, we only need to construct a function $\phi$ that belongs to $\mathcal{B}_1$ such that $\mathcal{J}(\phi)>0$. To this aim, define $\phi$ as
\begin{equation}\label{phi}
	\phi(x)=\left\{
	\begin{array}{lll}
		\frac{1}{\sqrt{h(O)\mu(O)}},&&x=O\\
		0,&&x\not= O.
	\end{array}\right.
\end{equation}
Obviously, we have $\phi\in\mathcal{H}$. Since $\int_Vh\phi^2d\mu=\mu(O)h(O)\phi^2(O)=1$, we also have  $\phi\in\mathcal{B}_1$. Furthermore, direct computations give us that
$$\frac{1}{p}\int_V|\phi|^pd\mu=\frac{1}{p}\mu(O)\phi^p(O)=\frac{1}{p}\frac{\mu(O)}{(h(O)\mu(O))^{p/2}}$$
and
\begin{eqnarray*}
	\frac{1}{2}\int_V|\nabla \phi|^2d\mu&=&\frac{1}{4}\sum_{x\in V}\sum_{y\sim x}w_{yx}(\phi(y)-\phi(x))^2\\
	&=&\frac{1}{4}\left\{\sum_{y\sim O}w_{yO}(\phi(y)-\phi(O))^2+\sum_{y\sim O}\sum_{z\sim y}w_{zy}(\phi(z)-\phi(y))^2\right\}\\
	&=&\frac{1}{4}\left\{{\text deg}(O)\phi^2(O)+{\text deg}(O)\phi^2(O)\right\}\\
	&=&\frac{1}{2}\frac{{\text deg}(O)}{h(O)\mu(O)}.
\end{eqnarray*}
Under the assumption (c3), one can easily check that
$$\mathcal{J}(\phi)=\frac{1}{p}\frac{\mu(O)}{(h(O)\mu(O))^{p/2}}-\frac{1}{2}\frac{{\text deg}(O)}{h(O)\mu(O)}>0.$$
Since $\phi\in\mathcal{B}_1$ and $\Lambda_1\geq \mathcal{J}(\phi)$, the lemma is proved.
$\hfill\Box$

\vspace{0.5cm}

\noindent {\bf Proof of Theorem \ref{existl}.}
Let us first deal with the special case $m=1$.
\vspace{0.5cm}

{\bf Case 1. $m=1$.}
Since $h_0=\inf_Vh>0$ and $\mu(x)\geq \mu_{\min}>0$, for any $u\in\mathcal{B}_1$ there holds
$$\|u\|_{\infty}\leq \frac{1}{\sqrt{\mu_{\min}h_0}}\left(\int_Vhu^2d\mu\right)^{1/2}\leq \frac{1}{\sqrt{\mu_{\min}h_0}}.$$
Hence we have
\begin{equation}\label{up}
	\int_V|u|^pd\mu\leq \|u\|_{\infty}^{p-2}\int_Vu^2d\mu\leq \frac{1}{h_0(\mu_{\min}h_0)^{(p-2)/2}}.
\end{equation}
This together with $\mathcal{J}(0)=0$ leads to
$$
0\leq \Lambda_1\leq\frac{1}{ph_0(\mu_{\min}h_0)^{(p-2)/2}}.$$
As a consequence, we can take a sequence of $\{u_n\}\subset\mathcal{B}_1$ (still denoted by $\{u_n\}$) such that $\mathcal{J}(u_n)\rightarrow \Lambda_1$ as $n\rightarrow\infty$.
This together with (\ref{up}) leads to
$$\int_V|\nabla u_n|^2d\mu=\frac{2}{p}\int_V|u_n|^pd\mu-2\Lambda_1+o_n(1)\leq C$$
for some constant $C$ depending only on $V$, $h$ and $p$. Then we know that $\{u_n\}$ is bounded in $\mathcal{H}$.
Since $\mu(x)\geq \mu_{\min}>0$,
$h\geq h_0>0$, and $h(x)\rightarrow+\infty$ as $\rho(x)\rightarrow+\infty$,
by Lemma 7 in \cite{LinYang2}, $\mathcal{H}$ is embedded in $L^p(V)$ compactly. Hence there exists a subsequence of $\{u_n\}$ (still denoted by $\{u_n\}$) and a function $u_*$ such that $u_k\rightarrow u_*$ in $L^p(V)$. For $u_*$, we have
$$\int_V|\nabla u_*|^2d\mu\leq\limsup_{n\rightarrow\infty}\int_V|\nabla u_n|^2d\mu$$
and
$$\int_Vhu_*^2d\mu\leq \limsup_{n\rightarrow\infty}\int_Vhu_n^2d\mu\leq 1.$$
Therefore, we can conclude that $u_*\in\mathcal{B}_1$ and it achieves the supremum of $\mathcal{J}$ in $\mathcal{B}_1$.

Since $|\nabla |u_*||\leq |\nabla u_*|$, we have $\mathcal{J}(|u_*|)\geq \mathcal{J}(u_*)$.
This immediately leads to $\mathcal{J}(|u_*|)=\Lambda_1$. Therefore, without loss of generality, we can assume $u_*\geq 0$ on $V$. Lemma \ref{bound} tells us $\mathcal{J}(u_*)=\Lambda_1>0$, which implies that $u_*\not\equiv 0$. We claim that
\begin{equation}\label{norm}
	\int_Vhu_*^2d\mu=1.
\end{equation}
If not, there must hold $\int_Vhu_*^2d\mu<1$. Let
$\widetilde{u}={(\int_Vhu_*^2d\mu)^{-1/2}}u_*$. It is easy to check $\widetilde{u}\in\mathcal{B}_1$ and
\begin{eqnarray*}
\mathcal{J}(\widetilde{u})&=&\frac{1}{p}\int_V|\widetilde{u}|^pd\mu-\frac{1}{2}\int_V|\nabla \widetilde{u}|^2d\mu\\
&=&\left(\int_Vhu_*^2d\mu\right)^{-p/2}\frac{1}{p}\int_V|{u}_*|^pd\mu-\left(\int_Vhu_*^2d\mu\right)^{-1}\frac{1}{2}\int_V|\nabla {u}_*|^2d\mu\\
&\geq&\left(\int_Vhu_*^2d\mu\right)^{-1}J(u_*)\\
&>&J(u_*)=\Lambda_1.
\end{eqnarray*}
This contradicts the definition of $\Lambda_1$. Then we know that (\ref{norm}) is true and $u_*\in \mathcal{S}_1\subset \mathcal{B}_1$. Therefore, $u_*$ is a maximizer of $\mathcal{J}$ under the constraint $u\in \mathcal{S}_1$. By a straightforward calculation, we get the Euler-Lagrange equation of $u_*$ as follows.
\begin{equation}\label{eqleq0}
	\left\{\begin{array}{lll}
	-\Delta u_*+\lambda hu_*=u_*^{p-1}&{\text in}& V\\
	u_*\geq 0&{\text in}& V\\
	\int_Vhu_*^2d\mu= 1,
    \end{array}\right.
\end{equation}
where $\lambda_1$ is an Euler-Lagrange multiplier written by
$$\lambda_1=\int_Vu_*^pd\mu-\int_V|\nabla u_*|^2d\mu,$$
which gives us that
\begin{eqnarray*}
\frac{\lambda_1}{2}&=&\frac{1}{2}\int_Vu_*^pd\mu-\frac{1}{2}\int_V|\nabla u_*|^2d\mu\\
&\geq&\frac{1}{p}\int_Vu_*^pd\mu-\frac{1}{2}\int_V|\nabla u_*|^2d\mu\\
&=&J(u_*)>0.
\end{eqnarray*}
Applying the maximum principle to (\ref{eqleq0}), as in the proof of Theorem \ref{existf}, we can confirm that $u_*(x)$ is strictly positive at any vertex $x\in V$.
This ends the proof for the case $m=1$. 
\vspace{0.5cm}

{\bf Case 2. $m>0$.}
If we make a transformation $u_m=\sqrt{m}v_m$ for $u_m\in \mathcal{H}$ and use $\mathcal{J}_m(v)$ to represent $\frac{m^{p/2}}{p}\int_V|v|^pd\mu-\frac{m}{2}\int_V|\nabla v|^2d\mu$, we have
\begin{eqnarray*}
	\Lambda_m&=&\sup_{u_m\in \mathcal{B}_m}\mathcal{J}(u_m)\\
	&=&\sup_{u_m\in\mathcal{B}_m}\left\{\frac{1}{p}\int_V|u_m|^pd\mu-\frac{1}{2}\int_V|\nabla u_m|^2d\mu\right\}\\
	&=&\sup_{v\in\mathcal{B}_1}\left\{\frac{m^{p/2}}{p}\int_V|v|^pd\mu-\frac{m}{2}\int_V|\nabla v|^2d\mu\right\}\\
	&=&\sup_{v\in \mathcal{B}_1}\mathcal{J}_m(v).
\end{eqnarray*}
Let $\phi$ be as in (\ref{phi}). Since $\phi\in \mathcal{B}_1$, in view of the condition (c3), we get $\Lambda_m>0$. Using a completely analogous argument as done
for the case $m=1$, we can find some $v_*\in\mathcal{S}_1$ such that $v_*\geq 0$ and
$$\mathcal{J}_m(v_*)=\Lambda_m=\sup_{v\in\mathcal{S}_1}\mathcal{J}_m(v).$$
Denote $\widetilde{u}=\sqrt{m}v_*$. There hold $\int_Vh\widetilde{u}^2d\mu=m$, $\widetilde{u}\geq 0$ and
$$\mathcal{J}(\widetilde{u})=\Lambda_m=\sup_{u\in\mathcal{S}_m}\mathcal{J}(u).$$
One can easily derive the Euler-Lagange equation of $\widetilde{u}$, which is
\begin{equation*}
	\left\{\begin{array}{lll}
	-\Delta \widetilde{u}+\lambda_m h\widetilde{u}=\widetilde{u}^{p-1}&{\text in}& V\\
	\widetilde{u}\geq0&{\text in}& V\\
	\int_Vh\widetilde{u}^2d\mu=m,
    \end{array}\right.
\end{equation*}
where
\begin{equation}\label{lambda}
	\lambda_m=\frac{1}{m}\left\{\int_V\widetilde{u}^pd\mu-\int_V|\nabla\widetilde{u}|^2d\mu\right\}\geq \frac{2}{m}\mathcal{J}(\widetilde{u})>0.
\end{equation}
By arguments similar to those in {\bf Case 1}, we can also confirm that $\widetilde{u}(x)>0$ for any $x\in V$ and this completes the proof of the theorem.
$\hfill\Box$
\vspace{0.5cm}

At the end of this section, we deal with the limit behaviour of the solution to \eqref{eql} on a locally finite graph as the mass $m$ varying.
\vspace{0.5cm}

\noindent{\bf Proof of Theorem \ref{limitl}.}
Suppose $u_m$ is a solution of \eqref{eql} with the mass constraint $m>0$. Let $v_m=u_m/\sqrt{m}$, then we have $\int_Vhv_m^2d\mu=1$. Since both $\mu$ and $h$ have positive lower bounds, there hold
$$\|v_m\|_{L^\infty(V)}\leq C \quad {\rm and} \quad \int_V|v_m|^sd\mu\leq C, \quad\forall s\geq 2,$$
for some positive constant $C$. By our assumption $h(x)\rightarrow+\infty$ as $\rho(x)\rightarrow+\infty$, up to a subsequence, we have
$$\lim_{m\rightarrow 0^+}\int_V(v_m-v)^2d\mu=0.$$
for some function $v$. In particular, $v_m(x)$ converges to $v(x)$ uniformly for $x\in V$ as $m\rightarrow 0^+$. Similarly, for the case $m\rightarrow +\infty$, we can also find some function $w$ such that $\lim_{m\rightarrow +\infty}\int_V(w_m-w)^2d\mu=0$ and $w_m(x)$ converges to $w(x)$ uniformly for $x\in V$ as $m\rightarrow +\infty$, where $w_m=u_m/\sqrt{m}$.
\vspace{0.5cm}

{\bf Case 1. $m\rightarrow 0^+$.}
If
$$\lambda_m=\frac{1}{m}\left(\int_V|u_m|^pd\mu-\int_V|\nabla u_m|^2d\mu\right)\geq 0,$$
we have
\begin{equation}\label{lim}
	\int_V|\nabla v_m|^2d\mu\leq m^{\frac{p}{2}-1}\int_Vv_m^pd\mu\leq Cm^{\frac{p}{2}-1}.
\end{equation}
This immediately leads to
$$\int_V|\nabla v|^2d\mu\leq \liminf_{m\rightarrow 0^+}\int_V|\nabla v_m|^2d\mu=0.$$
Hence $v\equiv c_0$ on $V$ for some constant $c_0$. If $c_0\not=0$, we have
$$+\infty=c_0^2\int_Vhd\mu=\int_Vhv^2d\mu\leq \liminf_{m\rightarrow 0^+}\int_Vhv_m^2d\mu=1,$$
which is a contradiction. Therefore $v\equiv 0$ on $V$.

On the other hand, if we have
$$\lambda_m=\frac{1}{m}\left(\int_V|u_m|^pd\mu-\int_V|\nabla u_m|^2d\mu\right)\geq -c,$$
for some positive constant $c$, instead of (\ref{lim}), we have
$$\int_V|\nabla v_m|^2d\mu\leq Cm^{\frac{p}{2}-1}+C,$$
where $C$ is some positive constant. This together with $\int_Vhv_m^2d\mu=1$ implies that $\{v_m\}$ is bounded in $\mathcal{H}$.
Since $\mathcal{H}$ is a reflexive Banach space, by Lemma 7 in \cite{LinYang2}, $\{v_m\}$ shall converge to some function $v$  weakly in $\mathcal{H}$, strongly in $L^s(V)$ for any $2\leq s\leq +\infty$,
as $m\rightarrow 0^+$. In view of \eqref{eql}, we have
\begin{equation}\label{v-eqn}
	\left\{\begin{array}{lll}
	-\Delta v_m+\lambda_m hv_m=m^{\frac{p}{2}-1}v_m^{p-1}&{\rm in}& V\\
	\int_Vhv_m^2d\mu=1.
    \end{array}\right.
\end{equation}
Since $\lambda_m=m^{\frac{p}{2}-1}\int_V|v_m|^{p-2}v_md\mu-\int_V|\nabla v_m|^2d\mu$ is bounded for $m\rightarrow 0^+$, we can assume that up to a subsequence, $-\lambda_m\rightarrow \lambda_0$ as $m\rightarrow 0^+$. As a consequence, $v\in \mathcal{H}$ satisfies
\begin{equation}\label{v-eq}
	\left\{\begin{array}{lll}
	-\Delta v=\lambda_0 hv&{\rm in}& V\\
	\int_Vhv^2d\mu\leq 1.
    \end{array}\right.
\end{equation}
Take $\{\phi_j\}\subset C_c(V)$ such that $\lim_{j\rightarrow \infty}\phi_j= v$ in $\mathcal{H}$. By testing the above equation by $\phi_j$, we obtain
$$\int_V\nabla v\nabla\phi_jd\mu=\lambda_0\int_Vhv\phi_jd\mu.$$
Let $j\rightarrow\infty$, we get
\begin{equation}\label{int-part}
	\int_V|\nabla v|^2d\mu=\lambda_0\int_Vhv^2d\mu.
\end{equation}
In view of (\ref{v-eq}), there holds either $v\equiv 0$, or $0<\int_Vhv^2d\mu\leq 1$. In the latter case, (\ref{int-part}) gives
$$\lambda_0=\frac{\int_V|\nabla v|^2d\mu}{\int_Vhv^2d\mu}>0.$$
\vspace{0.5cm}

{\bf Case 2. $m\rightarrow+\infty$.}
Since $w_m$ converges to $w$ point-wisely in $V$, we have
$$\int_{B_r(O)}(|\nabla w|^2+hw^2)d\mu\leq \liminf_{m\rightarrow+\infty}\int_V(|\nabla w_m|^2+hw_m^2)d\mu,$$
for any positive integer $r$, where $B_r(O)$ is the ball centred at $O\in V$ with the radius of $r$. Hence $w\in \mathcal{H}$. Noting that
$$\lambda_m=m^{\frac{p}{2}-1}\int_V|w_m|^pd\mu-\int_V|\nabla w_m|^2d\mu\geq -c,$$
we conclude
$$\lambda_m=O(m^{\frac{p}{2}-1})\quad{\text as}\quad m\rightarrow+\infty.$$
In view of (\ref{v-eqn}), one has
$$-m^{1-\frac{p}{2}}\Delta w_m=w_m^{p-1}-\frac{\lambda_m}{m^{\frac{p}{2}-1}} hw_m.$$
Let $m\rightarrow+\infty$, we arrive at
$$w^{p-1}=\lambda_\infty hw\quad{\rm in}\quad V,$$
where
$$\lambda_\infty=\lim_{m\rightarrow\infty}\frac{\lambda_m}{m^{\frac{p}{2}-1}}=\frac{\int_V|w|^pd\mu}{\int_Vhw^2d\mu}.$$
Thus the proof of the theorem is finished.
$\hfill\Box$

\begin{Rem}
	In the proof of Theorem \ref{existl}, we get the solution of \eqref{eql} by a maximization discussion and prove that the Lagrange multiplier $\lambda_m$ corresponding to the solution $u_m$ is positive (one can refer to \eqref{lambda}). On the other hand, there may exist solutions of \eqref{eql} with non-positive $\lambda_m$, which are not the same as the one in Theorem \ref{existl}. For this reason, in Theorem \ref{limitl}, we assume that $\lambda_m\geq -c$ for some constant $c\geq 0$ in order to deal with a more general case.
\end{Rem}

\section{Numerical results}

In this section, we illustrate our theoretical results by several numerical experiments. Since we can not really simulate a graph with infinite vertices or edges, we only carry on the experiments on a finite graph $G_6$ and verify the results in Theorem \ref{existf} and \ref{limitf}.

\begin{figure}[htbp]
	\centering
	\includegraphics[width=3.5in,height=2.5in]{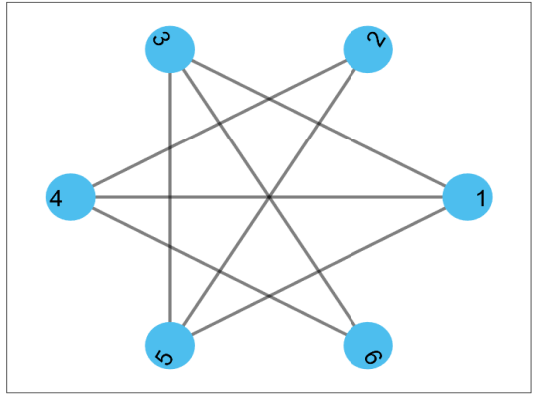}
	\caption{The graph $G_6$}
	\label{fig-g6}
\end{figure}

The graph $G_6$ has six vertices from $x_1$ to $x_6$ and its structure is shown in Figure \ref{fig-g6}. For simplicity, the symmetric weights of all edges of $G_6$ are set to be $1$. In order to give examples of the alternative results in Theorem \ref{existf} respectively, we adopt two settings for the positive measure $\mu(x)$, as shown in Table \ref{tab-mu}.

\begin{table}[!hbp]
	\centering
	\caption{The measure $\mu(x)$}
	\begin{tabular}{cccccc}
		\hline
		$\mu(x_1)$&$\mu(x_2)$&$\mu(x_3)$&$\mu(x_4)$&$\mu(x_5)$&$\mu(x_6)$ \\
		\hline
		\hline
		3   &  2   &  10   &  1   &  40   &  1     \\
		\hline
		$1$&$1$&$1$&$1$&$1$&$1$ \\
		\hline
	\end{tabular}
	\label{tab-mu}
\end{table}

We first compute the numerical solutions of the equation \eqref{eqf} for $p=3$ by MATLAB R2020a. The numerical solutions with the mass constraints $m=10^{-1}, 10,100$ and $\mu=(3,2,10,1,40,1)$ are shown in Table \ref{tab-sol}.

\begin{table}[!hbp]
	\centering
	\caption{The solutions with $m=10^{-1}, 10,100$ and $\mu=(3,2,10,1,40,1)$}
	\begin{tabular}{c|cccccc}
		\hline
		{\rm Mass}&$u(x_1)$&$u(x_2)$&$u(x_3)$&$u(x_4)$&$u(x_5)$&$u(x_6)$\\
		\hline
		\hline
		$10^{-1}$   &  0.0419  &  0.0419  &  0.0419  &  0.0419  &  0.0419  &0.0419  \\
		\hline
		$10$   &  0.1455  &  0.2252  &  0.0084  &  3.1068  &  0.0014  &0.4270  \\
		\hline
		$100$   &  0.1204  &  0.1817  &  0.0017  &  9.9881  &  0.0003  &0.3573  \\
		\hline
	\end{tabular}
	\label{tab-sol}
\end{table}

\begin{figure}[htbp]
	\centering
	\includegraphics[width=4.8in,height=2.4in]{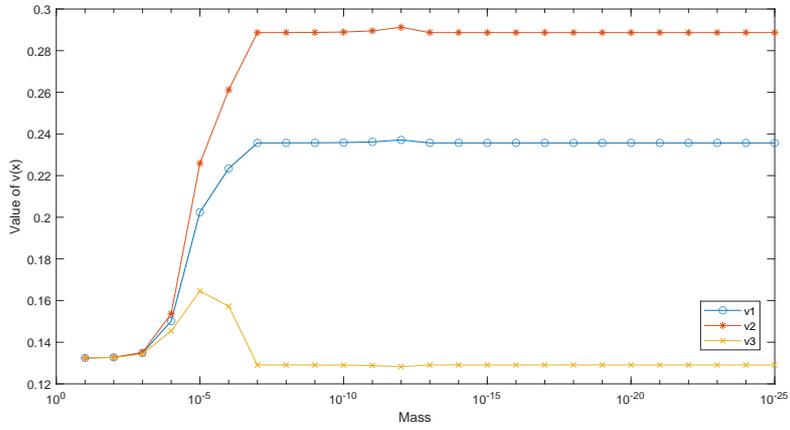}
	\caption{The values of $v_m(x)$ for $\mu=(3,2,10,1,40,1)$ as $m\rightarrow 0^+$}
	\label{fig-m0-sol}
\end{figure}

\begin{figure}[htbp]
	\centering
	\includegraphics[width=4.8in,height=2.4in]{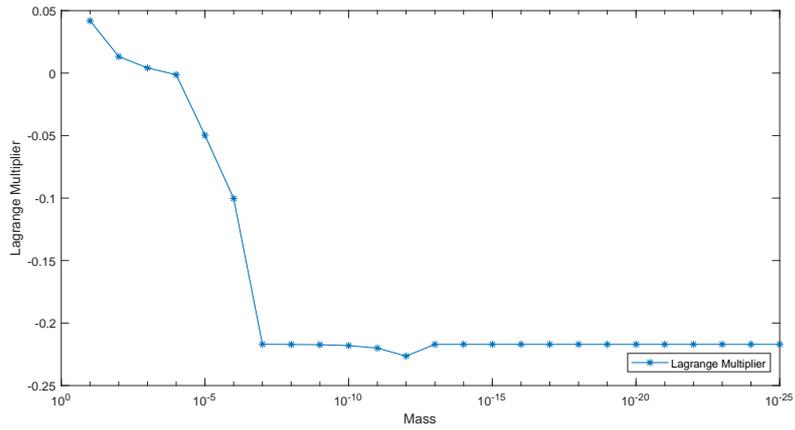}
	\caption{The Lagrange multiplier $\lambda_m$ for $\mu=(3,2,10,1,40,1)$ as $m\rightarrow 0^+$}
	\label{fig-m0-lm}
\end{figure}

Next, we simulate the situation $m\rightarrow 0^+$ in Theorem \ref{existf}. Let the mass constraint $m$ decrease from $10$ to $10^{-25}$, we find that for $\mu=(1,1,1,1,1,1)$, the rescaling solution $v_m=\frac{u_m}{\sqrt{m}}$ remains unchanged with $v_m\equiv|V|^{-\frac{1}{2}}=0.4082$. This at least gives us some hints about the conditions under which the first of the alternative result in Theorem \ref{existf} will occur. If we set $\mu=(3,2,10,1,40,1)$ and still let $m$ decrease from $10$ to $10^{-25}$, the other of the alternative result in Theorem \ref{existf} will occur. In order to show the curves of solutions more clearly, we only select the values of  $v_m$ at $x_1, x_2$ and $x_3$ to plot them in Figure \ref{fig-m0-sol}. The rescaling solution $v_m$ shall converge to a solution of the limit equation presented in Theorem \ref{existf}. To get the parameter $\lambda_0$ in the limit equation, we also compute and plot the Lagrange multiplier $\lambda_m$ as $m$ tends to $0^+$. We can find in Figure \ref{fig-m0-lm} that $\lambda_0=-\lim_{m\rightarrow 0^+} \lambda_m$ is a positive constant equals to $0.2171$, just as what Theorem \ref{existf} tells us. 

\begin{figure}[htbp]
	\centering
	\includegraphics[width=4.8in,height=2.4in]{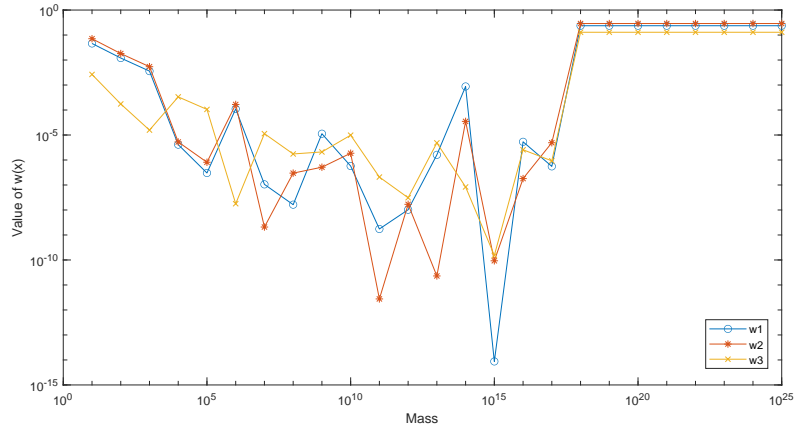}
	\caption{The values of $w_m(x)$ for $\mu=(3,2,10,1,40,1)$ as $m\rightarrow +\infty$}
	\label{fig-minfty-sol}
\end{figure}

\begin{figure}[htbp]
	\centering
	\includegraphics[width=4.8in,height=2.4in]{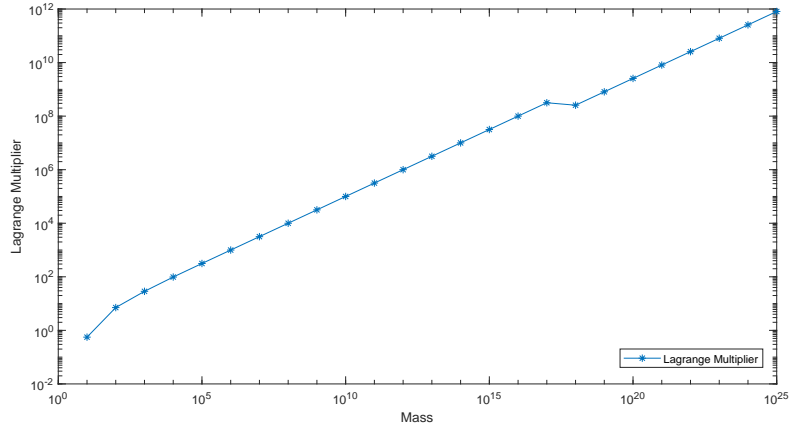}
	\caption{The Lagrange multiplier $\lambda_m$ for $\mu=(3,2,10,1,40,1)$ as $m\rightarrow +\infty$}
	\label{fig-minfty-lm}
\end{figure}

For the case $m\rightarrow +\infty$ in Theorem \ref{existf}, we let the mass constraint $m$ increase from $10$ to $10^{25}$. The values of the rescaling solution $w_m$ at $x_1, x_2$ and $x_3$ are plotted in Figure \ref{fig-minfty-sol}. For $m\geq 10^{18}$, the rescaling solution $w_m$ remains unchanged and satisfied the limit equation $w^{p-1}-\|w\|_p^pw=0$. Therefore, the values of $w_m$ for $m\geq 10^{18}$ are in fact the values of the limit function $w$ and we list them in Table \ref{tab-minfty-18}. The Lagrange multiplier $\lambda_m$ as $m$ tends to $+\infty$ is plotted in Figure \ref{fig-minfty-lm}. 

\begin{table}[!hbp]
	\centering
	\caption{The values of $w(x)$}
	\begin{tabular}{cccccc}
		\hline
		$w(x_1)$&$w(x_2)$&$w(x_3)$&$w(x_4)$&$w(x_5)$&$w(x_6)$ \\
		\hline
		\hline
		0.2357   &  0.2887   &  0.1291   &  0.4082   &  0.0645   &  0.4082     \\
		\hline
	\end{tabular}
	\label{tab-minfty-18}
\end{table}

According to the above results and their figures, we can confirm that the numerical experiments are completely consistent with what we have proved in Theorem \ref{existf} and \ref{limitf}.

\section*{Acknowledgements}
This research is supported by National Natural Science Foundation of China (No. 12271039).

\end{document}